\documentclass{amsart}
\usepackage{epsfig}
\newtheorem{lemma}{Lemma}
\newtheorem{thm}{Theorem}
\newtheorem{cor}{Corollary}
\newtheorem{prop}{Proposition}
\newtheorem{defn}{Definition}
\theoremstyle{definition}
\newtheorem{example}{Example}

\newtheorem{remark}{Remark}
\newcommand{\C}{{\mathbb C}}
\newcommand{\Z}{{\mathbb Z}}
\newcommand{\Q}{{\mathbb Q}}
\newcommand{\mbar}{\overline{m}}

\newcommand{\cpq}{\overline{\C P^2}}

\begin{document}

\title{Homology cobordism and classical knot invariants}
\author{Christian Bohr, Ronnie Lee}
\address{
Department of Mathematics, Yale University, P.O. Box 208283, New Haven, CT 06520-8283, 
}
\email{bohr@math.yale.edu}
\address{
Department of Mathematics, Yale University, P.O. Box 208283, New Haven, CT 06520-8283, 
}
\email{rlee@math.yale.edu}

\begin{abstract}
In this paper we define and investigate $\Z_2$--homology cobordism invariants of
$\Z_2$--homology 3--spheres which turn out to be related to classical invariants of knots.
As an application we show that many lens spaces have infinite
order in the $\Z_2$--homology cobordism group and we prove a lower bound
for the slice genus of a knot on which integral surgery yields a given $\Z_2$--homology sphere. 
We also give some
new examples of 3--manifolds which cannot be obtained by integral surgery on a knot.
\end{abstract}

\maketitle

\section{Introduction}

In recent years, gauge theoretical tools and new results in 4--dimensional topology 
have successfully been used to study the structure of the integral homology cobordism group. 
It is for instance a consequence of Donaldson's famous theorem about the intersection forms
of smooth 4--manifolds that the Poincar\'e homology sphere $\Sigma(2,3,5)$ has infinite order in this group, 
and M. Furuta found a family of Brieskorn spheres which generates a subgroup of infinite rank~\cite{Fu1}.
Recently N. Saveliev~\cite{Sa2} showed, using the $w$--invariant introduced in~\cite{FF}, that
a Brieskorn sphere with non--trivial Rokhlin invariant has infinite order in the
integral homology cobordism group.

However, many 3--manifolds arising naturally in knot theory, for instance double coverings of the 3--sphere
branched along knots, are not integral homology spheres, but still $\Z_2$--homology spheres. 
As in the case of integral homology spheres, the set of $\Z_2$--homology spheres
modulo the $\Z_2$--homology cobordism relation forms a group, the so called $\Z_2$--homology cobordism group
$\Theta^3_{\Z_2}$.
To study this group, we introduce, based on Furuta's result on the
intersection forms of smooth 4--dimensional spin manifolds~\cite{Fu2}, two invariants of 
$\Z_2$--homology spheres which turn out to be in fact invariants of the cobordism class
(see Theorem~\ref{basic}).
Exploiting that these invariants
are closely related to classical knot invariants like signature and slice genus, we prove estimates for them
in the case of lens spaces, which enables us to exhibit many examples of lens spaces which have infinite order in the
$\Z_2$--homology cobordism group, and we determine the slice genera of certain Montesinos knots.

A second relation between knot theory and cobordism classes of $\Z_2$--homology spheres is provided by the
simple fact that surgery along a knot with odd framing produces a $\Z_2$--homology sphere.
In this case our invariants can again be related to the slice genus of the knot. Using this,
we prove a lower bound for the slice genera of knots
on which integral surgery yields a given $\Z_2$--homology sphere. We also give 
new examples of 3--manifolds which cannot be obtained by integral surgery on a knot.

\section{Bordism invariants for $\Z_2$--homology spheres}

Recall that a closed connected and oriented 3--manifold $\Sigma$ is called a $\Z_2$--homology sphere if
$H_*(\Sigma;\Z_2)=H_*(S^3;\Z_2)$.
Two $\Z_2$--homology spheres $\Sigma_1,\Sigma_2$ are called $\Z_2$--homology cobordant  
if there exists a smooth 4--dimensional manifold $W$ with $\partial W=\Sigma_1 - \Sigma_2$
such that the inclusions $\Sigma_i \rightarrow W$ induce isomorphisms $H_*(\Sigma_i;\Z_2) \rightarrow H_*(W;\Z_2)$. 
The set of $\Z_2$--homology cobordism classes of $\Z_2$--homology spheres forms a group,
the so called $\Z_2$--homology cobordism group, which we denote by $\Theta^3_{\Z_2}$.
Addition in this group is given by taking the connected sum 
and the zero element is the equivalence class of the 3--sphere.

If $\Sigma$ is a $\Z_2$--homology sphere, there is a Rokhlin invariant $R(\Sigma) \in \Z_{16}$
which is defined to be the residue class of the signature of any (smooth) spin 4--manifold with boundary $\Sigma$.
It is easy to see that the Rokhlin invariant of a $\Z_2$--homology sphere is always even. 
The Rokhlin invariant defines a homomorphism
\[
R \colon \Theta_{\Z_2}^3 \rightarrow \Z_{16}
\]
with image $2 \Z_{16}$ (the Rokhlin invariant of $L(3,1)$ is two).

\begin{defn}
Let $\Sigma$ be a $\Z_2$--homology sphere. We define 
\begin{align*}
m(\Sigma) & = \max  \{ \frac{5}{4} \sigma(X) - b_2(X) \, \, \lvert \, \, w_2(X)=0, \partial X = \Sigma \}  \\
\mbar(\Sigma) & = \min  \{ \frac{5}{4} \sigma(X) + b_2(X) \, \, \lvert \, \, w_2(X)=0, \partial X = \Sigma\} 
\end{align*}
Here $X$ runs over all smooth spin 4--manifolds with boundary $\Sigma$.
\end{defn}

Surprisingly enough this simple definition actually yields bordism invariants which 
can be estimated if the manifold in question is a double covering of the 3--sphere branched along 
some knot. The following theorem summarizes some properties of the invariants $m$ and $\mbar$.

\begin{thm}\label{basic}
\mbox{}
\begin{enumerate}
\item\label{finite} 
The invariants $m$ and $\mbar$ are finite numbers. Moreover $m(\Sigma) \leq \mbar(\Sigma)$
where equality occurs if and only if
$m(\Sigma)=\mbar(\Sigma)=0$ and $R(\Sigma)=0$.
\item\label{orientation} For every $\Z_2$--homology sphere $\Sigma$, $m(-\Sigma)=-\mbar(\Sigma)$.
\item\label{sum}
If $\Sigma_1$ and $\Sigma_2$ are $\Z_2$--homology spheres, then 
\begin{align*}
m(\Sigma_1 \# \Sigma_2) \geq m(\Sigma_1) + m(\Sigma_2) \\
\mbar(\Sigma_1 \# \Sigma_2) \leq \mbar(\Sigma_1) + \mbar(\Sigma_2)
\end{align*}
\item\label{bordism}
If $\Sigma_1$ and $\Sigma_2$ are $\Z_2$--homology cobordant then we have
$m(\Sigma_1)=m(\Sigma_2)$ and $\mbar(\Sigma_1)=\mbar(\Sigma_2)$. In particular 
$m(\Sigma)=\mbar(\Sigma)=0$ if $\Sigma$ is the boundary of a $\Z_2$--acyclic 4--manifold.
\item\label{order} 
If $m(\Sigma) > 0$ or $\mbar(\Sigma) < 0$ then $\Sigma$ has infinite order in the
group $\Theta_{\Z_2}^3$. The same is true if $m(\Sigma)=0$ and $R(\Sigma) \neq 0$.
\item\label{knot}
If $\Sigma$ is a double covering of the 3--sphere branched along a knot $K$, then
\[
\frac{5}{4} \sigma(K) - 2 g^*(K) \leq m(\Sigma) \leq \mbar(\Sigma) \leq \frac{5}{4} \sigma(K) + 2g^*(K),
\]
where $g^*(K)$ denotes the slice genus of $K$ and $\sigma(K)$ its signature.
\end{enumerate}
\end{thm}

\begin{proof}
For the proof of statement~\ref{finite}, let $\Sigma$ be a $\Z_2$--homology sphere. 
Pick a spin manifold $X$ with boundary $\Sigma$. 
Now suppose we are given another spin manifold $Y$ such that $\partial Y=\Sigma$. Let $W=Y \cup (-X)$. 
As $\Sigma$ is a $\Z_2$--homology
sphere, $W$ is spin and $b_2(W)=b_2(X)+b_2(Y)$. Hence
\begin{align*}
\frac{5}{4} \sigma(Y) - b_2(Y) &= \frac{5}{4} (\sigma(Y)-\sigma(X)+\sigma(X)) - b_2(Y) \\
& = \frac{5}{4} (\sigma(Y)-\sigma(X)) + \frac{5}{4} \sigma(X) - b_2(Y) \\
& \leq b_2(W) - b_2(Y) + \frac{5}{4} \sigma(X) = b_2(X) + \frac{5}{4} \sigma(X)
\end{align*}
where we used Furuta's Theorem~\cite{Fu2} 
in the third line. Since $b_2(X) + \frac{5}{4} \sigma(X)$ does not depend on $Y$
we can conclude that $m(\Sigma)$ is finite and bounded from above by $\frac{5}{4} \sigma(X) + b_2(X)$.
Since this is true for every $X$ we also obtain that $\mbar(\Sigma) \geq m(\Sigma)$ and the finiteness of $\mbar$.
Now suppose we are given a $\Z_2$--homology sphere $\Sigma$ such that $m(\Sigma)=\mbar(\Sigma)$. Then
we can find spin 4--manifolds $X$ and $Y$ with $\partial X=\partial Y=\Sigma$
such that
\[
\mbar(\Sigma) = \frac{5}{4} \sigma(X) + b_2(X) =
\frac{5}{4} \sigma(Y) - b_2(Y) = m(\Sigma).
\]
Now consider the spin manifold $W=-X \cup Y$. Since $\sigma(W)=\sigma(Y)-\sigma(X)$ and
$b_2(W)=b_2(X)+b_2(Y)$, the above equality implies that $\frac{5}{4}\sigma(W)=b_2(W)$.
By Furuta's Theorem this is only possible if $b_2(W)=0$, and we obtain that $b_2(X)=b_2(Y)=0$.
Hence we have $m(\Sigma)=\mbar(\Sigma)=0$ and $R(\Sigma)=0$ as claimed.

Now let us proof property~\ref{orientation}. Pick a spin manifold $X$ with boundary $-\Sigma$ such that
$\frac{5}{4} \sigma(X) + b_2(X)=\mbar(-\Sigma)$. Then $-X$ bounds $\Sigma$, and we obtain that
\[
-\frac{5}{4} \sigma(X) - b_2(X) = - \mbar(-\Sigma) \leq m(\Sigma).
\]
Hence we have $-m(\Sigma) \leq \mbar(-\Sigma)$.
If we choose a spin manifold $Y$ with boundary $\Sigma$ such that $\frac{5}{4} \sigma(Y) - b_2(Y)=m(\Sigma)$,
then $-Y$ bounds $-\Sigma$ and we have
\[
-\frac{5}{4} \sigma(Y) + b_2(Y) = - m(\Sigma) \geq \mbar(-\Sigma).
\]
These two inequalities imply the desired result.

As to assertion~\ref{sum} first note that -- thanks to property ~\ref{orientation} -- it suffices to prove this for
$m$. Choose spin manifolds $X_1$ and $X_2$ with boundaries $\Sigma_1$ and $\Sigma_2$ such that
$\frac{5}{4} \sigma(X_i) - b_2(X_i) = m(\Sigma_i)$. Then $V=X_1 \natural X_2$ 
is spin with boundary $\Sigma_1 \# \Sigma_2$
and $\frac{5}{4} \sigma(V) - b_2(V) = m(\Sigma_1) + m(\Sigma_2)$.

To prove assertion~\ref{bordism} note that again we only have to prove the required property for $m$.
The condition $[\Sigma_1]=[\Sigma_2]$ implies that there exists
a spin 4--manifold $W$ with boundary $\partial W=\Sigma_1 - \Sigma_2$ such that $\sigma(W)=b_2(W)=0$. 
Pick manifolds $X_1,X_2$
such that $\partial X_i=\Sigma_i$ and $\frac{5}{4} \sigma(X_i)-b_2(X_i)=m(\Sigma_i)$. Consider the
spin manifolds 
\begin{align*}
V_1 &=X_1 \cup _{\Sigma_1} (- W) \\
V_2 &=X_2 \cup _{\Sigma_2} W 
\end{align*}
Then $\partial V_1=\Sigma_2$,
$\sigma(V_1)=\sigma(X_1)$ and $b_2(V_1)=b_2(X_1)$. Similarly $\partial V_2=\Sigma_1$, $\sigma(V_2)=\sigma(X_2)$ and
$b_2(V_2)=b_2(X_2)$. So we obtain
\begin{align*}
m(\Sigma_2) \geq \frac{5}{4} \sigma(V_1) - b_2(V_1) = m(\Sigma_1) \\
m(\Sigma_1) \geq \frac{5}{4} \sigma(V_2) -b_2(V_2) = m(\Sigma_2)
\end{align*}
and the claim follows. To prove the second part of the assertion we only have to prove that
$m(S^3)=\mbar(S^3)=0$. Since $b_2(D^4)=\sigma(D_4)=0$ we have $m(S^3) \geq 0$ and $\mbar(S^3) \leq 0$. 
By statement~\ref{finite}, we obtain 
\[
0 \leq m(S^3) \leq \mbar(S^3) \leq 0
\]
and therefore $m(S^3)=\mbar(S^3)=0$.

The first part of statement~\ref{order} is an immediate consequence of \ref{bordism} and
\ref{sum}. As to the second part, the assumption $m(\Sigma)=0$ implies that there exists a spin manifold $X$
with boundary $\Sigma$
such that $\frac{5}{4}\sigma(X)=b_2(X)$. Now suppose there is some $n >0$ such that $n\Sigma$ is the boundary
of a $\Z_2$--acyclic manifold $V$. Consider the spin manifold
\[
Y=\underbrace{X \cup X \cdots X}_{n \text{ copies }} \cup -V.
\]
Then $\frac{5}{4}\sigma(Y)=b_2(Y)$ which, by Furuta's Theorem, implies that $\sigma(Y)=0$. But $\sigma(X)=n\sigma(X)$,
hence we obtain that $\sigma(X)=0$, in contradiction to $R(\Sigma) \neq 0$.

To prove property~\ref{knot}, note that given a surface $F$ in $D^4$ with boundary $K$, there is a 
double covering $X \rightarrow D^4$ branched along $F$ with boundary $\Sigma$. The manifold $X$ is spin,
and it is well known that $b_2(X)=2g(F)$ whereas $\sigma(X)=\sigma(K)$ (see~\cite{CL} for a proof). 
This implies the claimed inequalities.
\end{proof}

\begin{remark} Observe that the invariants $m$ and $\mbar$ are both not additive. 
In fact, if one of them were additive,
then property~\ref{orientation} would imply that $m(\Sigma)=\mbar(\Sigma)$ for all $\Z_2$--homology spheres $\Sigma$.
As there are clearly $\Z_2$--homology spheres with non--zero Rokhlin invariant -- for instance $\Sigma=L(3,1)$ --
this is impossible by property~\ref{finite}.

We also note that the invariants $m$ and $\mbar$ are in general not integers (however they are always multiples of 
$\frac{1}{2}$) and that they are related by the formula
\[
m(\Sigma) = \mbar(\Sigma) = \frac{1}{4} R(\Sigma) \mod 2\Z
\]
to the Rokhlin invariant. In particular $m(\Sigma)-\mbar(\Sigma)$ is always an even integer.
\end{remark}

\begin{example}\label{mpqr}
Suppose we are given positive integers $p \leq q \leq r$ such that exactly one of these three numbers is even,
$\frac{1}{p} + \frac{1}{q} + \frac{1}{r} < 1$ and $p+q+r \leq 22$.
Let $T_{p,q,r}$ denote the weighted graph shown in figure~\ref{tpqr} and consider the 4--manifold $X_{p,q,r}$ 
obtained by plumbing according to $T_{p,q,r}$. It is not hard to check that the determinant of $T_{p,q,r}$ is
-- up to sign -- the number $pqr-pq-pr-qr$ and is therefore odd, hence the boundary 
$\Sigma_{p,q,r}=\partial X_{p,q,r}$ is a $\Z_2$--homology sphere.

\begin{figure}[ht]
\begin{center}
\epsfig{file=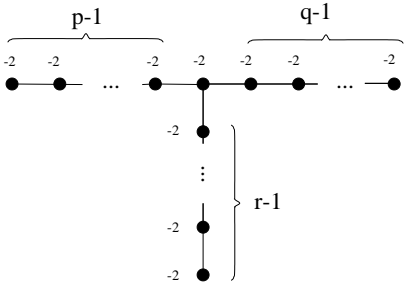}
\caption{The graph $T_{p,q,r}$\label{tpqr}}
\end{center}
\end{figure}

As the signature of $T_{p,q,r}$ is $4-p-q-r$ and the rank is $p+q+r-2$, we obtain that
\[
\mbar(\Sigma_{p,q,r}) \leq -\frac{1}{4} (p+q+r) + 3.
\]
Now it has been proved in~\cite{L} and~\cite{P} that the bilinear form $T_{p,q,r}$ can be realized by a collection
of $(-2)$--curves in a K3--surface, and a regular neighborhood of such a configuration of spheres is diffeomorphic
to $X_{p,q,r}$, hence we have an embedding $X_{p,q,r} \subset K3$. Let $Y$ denote the closure of the complement. 
Then $\partial Y =-\Sigma_{p,q,r}$, $sign(Y)=-16-sign(X_{p,q,r})=p+q+r-20$ and 
$b_2(Y)=22-b_2(X_{p,q,r})=24-p-q-r$, and we obtain
\[
m(\Sigma_{p,q,r}) \geq -\frac{5}{4}(p+q+r-20)-24+p+q+r = -\frac{1}{4}(p+q+r)+1.
\]

Now we claim that $\mbar(\Sigma_{p,q,r}) > m(\Sigma_{p,q,r})$. In fact, if $p+q+r \not\equiv 4 \mod 16$, 
this follows from
Theorem~\ref{basic}, as in this case $R(\Sigma_{p,q,r}) \neq 0$. Our conditions on $p,q,r$ exclude the case that
$p+q+r=4$, so the only remaining case we have to check is $p+q+r=20$. But then we already know that 
$\mbar(\Sigma_{p,q,r}) \leq -2$, and therefore we can again conclude that 
$\mbar(\Sigma_{p,q,r}) > m(\Sigma_{p,q,r})$, which proves our claim.

As the difference $m-\mbar$ is always an even integer, this discussion shows that the 
above estimates for $m(\Sigma_{p,q,r})$
and $\mbar(\Sigma_{p,q,r})$ are sharp, i.e.
\[
m(\Sigma_{p,q,r})=-\frac{1}{4}(p+q+r)+1=\mbar(\Sigma_{p,q,r})-2.
\]
Note that $\Sigma_{2,3,7}$ is minus the Brieskorn sphere $\Sigma(2,3,7)$, hence we can conclude that 
$m(\Sigma(2,3,7))=0$. As $R(\Sigma(2,3,7))=8 \neq 0$, we
obtain that $\Sigma(2,3,7)$ has infinite order in the group $\Theta^3_{\Z_2}$ 
although it bounds a rational ball~\cite{FS2}. This also follows from the results in~\cite{Sa3}.
\end{example}

As we have seen in Example~\ref{mpqr}, the Brieskorn sphere $\Sigma(2,3,7)$ is an element of infinite order
in the kernel of the natural map 
$\Theta^3_{\Z_2} \rightarrow \Theta^3_{\Q}$, where the latter group is defined in the obvious way. Also note that
this map is not onto, as the Rokhlin invariant of every $\Z_2$--homology sphere is even. A similar argument shows
that $\Theta^3_{\Z} \rightarrow \Theta^3_{\Z_2}$ is not onto. It is conceivable that this map is also not 
one--to--one, and it would be interesting to find explicit examples of integral homology spheres which are not the
boundary of a $\Z$--acyclic 4--manifold but are bounded by some $\Z_2$--acyclic manifold.
We remark that
Furuta's arguments given in~\cite{Fu1} actually show that the family $\Sigma(2k+1,4k+1,4k+3)$ of Brieskorn spheres 
generates a subgroup of infinite order in $\Theta^3_{\Z_2}$ (this has been used by H. Endo~\cite{E} to prove
that the corresponding family of Pretzel knots spans an infinite dimensional subspace in the smooth knot
concordance group), so the restriction of $\Theta^3_{\Z} \rightarrow \Theta^3_{\Z_2}$ to this subgroup is
one--to--one.

\begin{cor}\label{montesinos}
Suppose we are given positive integers $p \leq q \leq r$ such that exactly one of these three numbers is even,
$\frac{1}{p} + \frac{1}{q} + \frac{1}{r} < 1$ and $p+q+r \leq 22$.
Then the Montesinos knot ${\mathfrak m}(2;(p,p-1),(q,q-1),(r,r-1))$ has slice genus
and unknotting number $\frac{1}{2}(p+q+r)-1$.
\end{cor}

\begin{proof}
Let $K$ denote the knot ${\mathfrak m}(2;(p,p-1),(q,q-1),(r,r-1))$.
The double covering $\Sigma$ of $S^3$ branched along this knot is $-\Sigma_{p,q,r}$, hence 
\[
m(\Sigma)=-\mbar(\Sigma_{p,q,r})=\frac{1}{4}(p+q+r)-3
\] 
by Example~\ref{mpqr}.
It is not hard to verify that $\sigma(K)=p+q+r-4$. By Theorem~\ref{basic}
we obtain that
\[
2g^*(K) \geq \frac{5}{4} \sigma(K) - m(\Sigma)= \frac{5}{4}(p+q+r)-5-\frac{1}{4}(p+q+r)+3=p+q+r-2,
\]
so $g^*(K) \geq \frac{1}{2}(p+q+r)-1$. Now the knot $K$ is build up from the rational tangles
$t(p,p-1),t(q,q-1)$ and $t(r,r-1)$. 
The rational tangle $t(p,p-1)$ can be changed to the trivial tangle by $\frac{p}{2}$ crossing changes 
if $p$ is even and by $\frac{p-1}{2}$ crossing changes if $p$ is odd, so we obtain that 
$u(K) \leq \frac{1}{2}(p+q+r)-1$.
Of course $g^*(K) \leq u(K)$, and therefore we can conclude that $u(K)=g^*(K)=\frac{1}{2}(p+q+r)-1$.
\end{proof}

Note that the Montesinos knots in Corollary~\ref{montesinos} are examples of knots where the
inequality $g^*(K) \geq \frac{1}{2} |\sigma(K)|$ from~\cite{M}
is not sharp. In some sense the last statement of Theorem~\ref{basic} can be seen as a 
refinement of this inequality involving the cobordism class of the double branched covering.

\begin{example} Suppose that $\Sigma$ is an integral homology sphere with non--zero Rokhlin invariant 
which is the result of rational surgery on a torus knot. Then $m(\Sigma) \geq 0$. 
This follows from~\cite{Sa3} where Saveliev shows 
that certain Brieskorn spheres bound spin 4--manifolds $X$ with
intersection forms $-aE_8+bH$,$ a \geq b$ and uses this to prove that they have infinite order in the homology
cobordism group. For these homology spheres $m(\Sigma) \geq \frac{5}{4} \sigma(X) - b_2(X)=2(a-b) \geq 0$.
\end{example}
 
\section{Some computations for lens spaces}

As A. Casson and J. Harer demonstrated in~\cite{CH}, all lens spaces of the form $L(t^2,qt+1)$ with coprime
numbers $t$ and $q$ where $t$ is odd
bound $\Z_2$--acyclic 4--manifolds, hence they are zero in the group $\Theta^3_{\Z_2}$.
In this section, we prove estimates for the invariants $m$ and $\mbar$ of lens spaces
and use them to find many examples of lens spaces which have in fact infinite order in the $\Z_2$--homology
cobordism group.
In particular we will obtain a complete list of those lens spaces $L(\alpha,\beta)$ with odd $\alpha \leq 11$
whose cobordism classes have infinite order.

It is well known that a lens space $L(\alpha,\beta)$ is a double covering of the
3--sphere branched along the two--bridge link $S(\alpha,\beta)$, see for instance~\cite{BZ}. To be able to
apply statement~\ref{knot} in Theorem~\ref{basic}, we have to compute respectively to estimate the signature
and the slice genus of such a link.
Let us start by fixing some notations.

\begin{figure}[ht]
\begin{center}
\epsfig{file=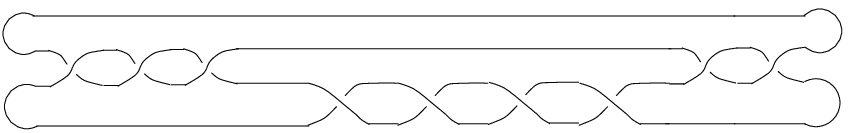}
\caption{P(3,2,2)\label{twobridgeknot}}
\end{center}
\end{figure}

\begin{defn} Assume that we are given a sequence $(a_1,b_1,a_2, \cdots, a_n)$ of $2n-1$ non--zero integers.
Then we define the link $P(a_1,b_1, \cdots, a_n)$ to be the 4--plat
obtained by closing the 3--string braid $\sigma_2^{-a_1}\sigma_1^{2b_1} \cdots \sigma_2^{-a_n}$ as pictured
in figure~\ref{twobridgeknot}
(here we use the convention that the generators of the braid group have positive crossings if the
two strings involved have parallel orientations).
\end{defn}

The link $P(a_1,b_1, \cdots, a_n)$ is a knot if and only if $\sum_i a_i$ is odd, otherwise it
is a two component link.

\begin{defn}
Assume that we are given coprime integers $0 < \beta < \alpha$ such that $\beta$ is odd. Then an admissible
continuous fractions decomposition is a presentation of $\frac{\alpha}{\beta}$ as a continuous fraction
\[
\frac{\alpha}{\beta} =[a_1,2b_1,a_2, \cdots a_n] = 
a_1 + \cfrac{1}{2b_1 + \cfrac{1}{a_2+ \cfrac{1}{2b_2+ \cdots }}}
\]
where the $a_i$ and $b_i$ are integers such that $a_i b_i > 0$ for all $i=1, \cdots, n-1$.
\end{defn}

Observe that such an admissible continuous fractions decomposition exists 
for every pair $0 < \beta < \alpha$ with $gcd(\alpha,\beta)=1$ and $\beta \equiv 1 \mod 2$.
It is also well known~\cite{BZ} that the 4--plat $P(a_1,b_1, \cdots, a_n)$ is nothing else than the
2--bridge link $S(\alpha,\beta)$. This description of 2--bridge links turns out to be particularly useful for
computing the signature and slice genus of such a link.
The following fact can be found in \cite{B}.

\begin{lemma}[see \cite{B}]\label{signaturecomputation}
Assume that $0 < \beta < \alpha$ are coprime integers and that $\beta$ is odd. Pick an admissible 
continuous fractions expansion
\[
\frac{\alpha}{\beta} = [a_1,2b_1, \cdots, a_n].
\]
Then the signature of the two--bridge link $S(\alpha,\beta)$ is given by
\[
\sigma(S(\alpha,\beta))= \sum_i a_i - \frac{a_n}{|a_n|}.
\]
\end{lemma}

To find an upper bound for the slice genus of a two--bridge knot we will use the following general observation.

\begin{prop}\label{slicegenus}
Suppose that $K_1$ and $K_2$ are links such that $K_2$ is obtained from $K_1$ by $p$ positive and $n$
negative crossing changes. Then
\[
g^*(K_1) \leq g^*(K_2) + \max \{p,n\}.
\]
\end{prop}

\begin{proof}
The trace of a homotopy given by the crossing changes is a union of immersed annuli $A$ in $S^3 \times I$
with $p$ positive and $n$ negative self--intersection points
connecting $K_1 \subset S^3 \times \{0\}$ and $K_2 \subset S^3 \times \{1\}$
(see \cite{CL}). Pick a connected surface $F \subset D^4$
with boundary $K_2$ such that $g(F)=g^*(K_2)$. By gluing this surface with $A$ we obtain a connected immersed
surface $F' \subset D^4=(S^3 \times I) \cup D^4$ having boundary $K_1$ whose genus is $g(F)$ and which has $p$
positive and $n$ negative self intersection points. Now let us assume that $p \geq n$. Since we can join two
self--intersection points of opposite signs by a handle we can construct a surface bounding $K_1$
which has genus $g(F)+n$ and $p-n$ positive self--intersection points. Replacing the remaining self--intersections
points by $p-n$ handles, we end up with an embedded surface with boundary $K_1$ which 
has genus $g(F)+p$ and the claim follows. In the case that $n \geq p$ a similar argument applies.
\end{proof}

\begin{lemma}\label{genuscomputation}
Suppose that we have a sequence $(a_1,b_1,a_2, \cdots, a_n)$ of non--zero integers
such that $\sum_i a_i$ is odd.
Define numbers $o^+,o^-$ by
\begin{align*}
o^+ &= \# \{1 \leq i \leq n \, \, \lvert \, \, a_i \equiv 1 \!\! \mod 2, a_i > 0 \}, \\
o^- &= \# \{1 \leq i \leq n \, \, \lvert \, \, a_i \equiv 1 \!\! \mod 2, a_i < 0 \}.
\end{align*}
Then the slice genus
$g^*$ of the knot $P(a_1, b_1, \cdots, a_n)$ is bounded from above by
\[
g^* \leq 
\frac{1}{4} \max \{ \ \sum_i (|a_i|-a_i) + 2 o^+ - 2, 
\sum_i (|a_i|+a_i) + 2 o^- - 2 \}.
\]
\end{lemma}

\begin{proof}
Suppose that we have an index
$0 \leq i \leq n$ such that $a_i>0$. If $a_i$ is odd, we can
deform the knot $P(a_1,b_1, \cdots, a_n)$ into the knot defined by the sequence
$(a_1,b_1, \cdots, b_{i-1},1,b_i,a_{i+1}, \cdots, a_n)$ by performing $\frac{1}{2}(a_i-1)$ negative crossing changes. 
If $a_i>0$ is even, we can do $\frac{1}{2}a_i$ negative crossing changes to obtain 
$P(a_1,b_1, \cdots, b_{i-1},0,b_i,a_{i+1}, \cdots, a_n)$. Repeating this for every index $i$ for which $a_i$ is
positive, we eventually obtain a knot $P(a_1',b_1,a_2',b_2, \cdots, a_n)$ for which $a_i'=a_i$ if $a_i < 0$
and $a_i' \in \{0,1\}$ otherwise after having performed
\[
n=\sum_{\substack{i \\ a_i >0 \\ a_i \text{ odd}}} \frac{1}{2}(a_i-1) + 
\sum_{\substack{i \\ a_i >0 \\ a_i \text{ even}}} \frac{1}{2} a_i
= \frac{1}{4} (\sum_i (|a_i|+a_i) - 2 o^+)
\]
negative crossing changes. 

A similar reduction can be done if $a_i<0$. In this case we can do $\frac{1}{2}(|a_i|-1)$ respectively
$\frac{1}{2}|a_i|$ positive crossing changes, depending on whether $a_i$ is even or odd. So we see that 
after performing 
\[
p=\frac{1}{4} (\sum_i (|a_i|-a_i) - 2 o^-)
\]
additional positive crossing changes, we end up with a link $P(a_1'',b_1,a_2'',b_2, a_n'')$ where
$a_i'' \in \{-1,0,1\}$ and $a_i''=0$ if and only if $a_i$ is even. 

Observe that the knot $P(a_1'',b_1,a_2'',b_2, a_n'')$
is the boundary of an obvious Seifert surface
which has genus $\frac{1}{2} (o^+ + o^- -1)$. 

Now assume that $p \geq n$, i.e. $p=\max \{p,n \}$. By Proposition~\ref{slicegenus} we can conclude that
\begin{align*}
g^*(S(\alpha,\beta))  \leq p + g(F) & =\frac{1}{4} (\sum_i (|a_i|-a_i) - 2 o^-) + \frac{1}{2} (o^+ + o^- - 1) \\
& = \frac{1}{4} (\sum_i (|a_i|-a_i) + 2 o^+ - 2). 
\end{align*}
If we have $n \geq p$, we can use the same argument to obtain the lower bound
\begin{align*}
g^*(S(\alpha,\beta))  \leq n &+ g(F) =\frac{1}{4} (\sum_i (|a_i|+a_i) - 2 o^+) + \frac{1}{2} (o^+ + o^- - 1) \\
& = \frac{1}{4} (\sum_i (|a_i|+a_i) + 2 o^- - 2).
\end{align*}
Since of course either $p \geq n$ or $n \geq p$ the claimed inequality follows.
\end{proof}

\begin{prop}\label{lens}
Suppose that $0 < \beta < \alpha$ are coprime odd numbers and that we are given a continuous fractions decomposition
\[
\frac{\alpha}{\beta} =[a_1,2b_1,a_2, \cdots a_n] 
\]
with non--zero integers $a_i, b_i$ such that $a_i b_i > 0$ for all $i$. Let
\begin{align*}
o^+ &= \# \{1 \leq i \leq n \, \, \lvert \, \, a_i \equiv 1 \!\! \mod 2, a_i > 0 \}, \\
o^- &= \# \{1 \leq i \leq n \, \, \lvert \, \, a_i \equiv 1 \!\! \mod 2, a_i < 0 \}.
\end{align*}
Then
\begin{align*}
m(L(\alpha,\beta))  \geq 
\frac{5}{4} (\sum_{i=1}^n a_i - \frac{a_n}{|a_n|})
- \frac{1}{2} \max \Big\{ & \sum_{i=1}^n (|a_i|+a_i)+2o^- -2, \\
& \sum_{i=1}^n (|a_i|-a_i)+2o^+ -2  \Big\} .
\end{align*}
and
\begin{align*}
\mbar(L(\alpha,\beta))  \leq 
\frac{5}{4} ( \sum_{i=1}^n a_i - \frac{a_n}{|a_n|} )
+ \frac{1}{2} \max \Big\{ & \sum_{i=1}^n (|a_i|+a_i)+2o^- -2, \\
& \sum_{i=1}^n (|a_i|-a_i)+2o^+ -2  \Big\} .
\end{align*}
\end{prop}

\begin{proof}
It is well known that the lens space $L(\alpha,\beta)$ is a double covering of $S^3$, branched along the 
two--bridge knot $S(\alpha,\beta)=P(a_1, b_1, \cdots, a_n)$, see for instance \cite{BZ}. 
By Lemma~\ref{genuscomputation}, there exists a surface
$F \subset D^4$ with boundary $S(\alpha,\beta)$ which has genus
\[
g(F)=
\frac{1}{4} \max \{ \ \sum_i (|a_i|-a_i) + 2 o^+ - 2 , 
\sum_i (|a_i|+a_i + 2 o^- - 2) \}.
\]
Using the expression for the
signature of $S(\alpha,\beta)$ derived in Lemma~\ref{signaturecomputation}, 
Theorem~\ref{basic} now immediately yields the desired result.
\end{proof}

\begin{remark}
As to the orientation of lens spaces, we are using the convention from \cite{BZ} that the oriented lens space
$L(\alpha,\beta)$ for $0 < \beta < \alpha$ is the double covering of the 3--sphere branched along the 
two--bridge link $S(\alpha,\beta)$. With this choice of orientations, 
the covering of the left--handed trefoil knot is $L(3,1)$. 
Note that this is in accordance with the convention used in \cite{GS}, where the lens space $L(\alpha,\beta)$
is defined to be the result of rational surgery along the unknot with framing $\frac{-\alpha}{\beta}$
(see \cite{GS}, Exercise 6.3.5).
\end{remark}

By computing the terms appearing in the statement of Proposition~\ref{lens} 
we can find many examples of lens spaces which have infinite order
in the $\Z_2$--homology cobordism group. The following table shows -- up to orientation -- all the lens
spaces whose first homology groups have odd order less or equal than 13 together with the estimates for 
$m$ and $\mbar$ provided by Proposition~\ref{lens},  the continuous fractions decomposition used
for the computation and the order of the lens space in $\Theta^3_{\Z_2}$ as far as it is known (note that the lens
spaces $L(5,3)$ and $L(13,5)$ have orientation reversing diffeomorphisms, the fact that $L(9,5)$ is the boundary
of a $\Z_2$--acyclic manifold is proved in~\cite{CH})).

\begin{table}
\begin{center}
\begin{tabular}{|r|r|r|r|r|}
\hline
$L(\alpha,\beta)$ & $m \geq$ & $\mbar \leq$ & $[a_1,2b_2,a_2, \cdots ,a_n]$  & order \\
\hline
(3,1) & 0.5 & 4.5 & [3] & $\infty$ \\
(5,3) & -2.0 & 2.0 & [1,2] & $\leq 2$ \\
(7,1) & 1.5 & 13.5 & [7] & $\infty$ \\
(7,3) & 0.5 & 4.5 & [2,4,-1] & $\infty$ \\
(9,1) & 2.0 & 18.0 & [9] & $\infty$ \\
(9,5) & -2.0 & 2.0 & [1,2,-1,-2,-1] & 0 \\
(11,1) & 2.5 & 22.5 & [11] & $\infty$ \\
(11,3) & 0.5 & 4.5 & [3,2,-2] & $\infty$ \\
(11,5) & 0.5 & 4.5 & [2,6,-1] & $\infty$ \\
(13,1) & 3.0 & 27.0 & [13] & $\infty$ \\
(13,3) & 1.0 & 9.0 & [4,4,1] & $\infty$ \\
(13,5) & -2.0 & 2.0 & [2,2,-3] & $\leq 2$\\
(13,7) & -2.0 & 2.0 & [1,2,-1,-4,-1] & ? \\
\hline
\end{tabular}
\end{center}
\caption{\label{table1}Bounds for $m(L(p,q))$}
\end{table}

\begin{cor}\label{positive}
Suppose that $0 < \beta < \alpha$ are coprime odd numbers and that $\frac{\alpha}{\beta}$ has a continuous fractions
decomposition
\[
\frac{\alpha}{\beta}=[a_1,2b_1,a_2,2b_2,a_3, \cdots,a_n]
\]
such that $a_i, b_i >0$ for all $i$. Then the lens space $L(\alpha,\beta)$ has infinite order
in the $\Z_2$--homology cobordism group. 
\end{cor}

\begin{example}
Suppose $n>0$ is some positive integer. Then the lens space $L(10n+1,8n+1)$
has infinite order in $\Theta^3_{\Z_2}$. In fact, first note that $10n+1$ and
$8n+1$ are odd and coprime (if some prime $p$ divides $10n+1$ and $8n+1$, then $2a \equiv 0 \mod p$, 
in contradiction to $10n+1 \equiv 0 \mod p$). A continuous fractions decomposition is given by
\[
\frac{10n+1}{8n+1} = 1 + \frac{2n}{8n+1} = 1 + \frac{1}{4+\frac{1}{2n}}=[1,4,2n].
\]
Since all the coefficients are positive Corollary~\ref{positive} applies and the claim 
follows. 

As the lens spaces $L(10n+1,8n+1)$ and $L(10m+1,8m+1)$ are diffeomorphic if and only if
$n=m$, we obtain in particular that there
is an infinite family of lens spaces which have infinite order in the $\Z_2$--homology cobordism group.
\end{example}

\begin{example}
We have seen that a $\Z_2$--homology sphere obtained by integral surgery on the unknot -- namely a lens space
$L(n,1)$ -- has infinite order in the $\Z_2$--homology cobordism group unless it is an integral homology sphere,
in which case it is $\Z_2$--homology cobordant to $S^3$.
It turns out that the same is true if the knot is slice, unless the framing is $\pm 1$. In fact,
suppose that $\Sigma$ is obtained by doing surgery with odd framing $n$ on a slice knot $K$. Without loss
of generality we can assume that $n > 0$. Let $W$ be the trace of the surgery, i.e. $W$ is obtained from $D^4$
by attaching a 2--handle along $K$ with framing $n$. Then an embedded disk in the 4--ball with boundary $K$
and the core of the 2--handle can be glued together to given an embedded sphere $S \subset W$ with self--intersection
number $n$ which generates $H_2(W;\Z)$. 
Let $X$ be the manifold obtained from $W$ by removing a tubular neighborhood of $S$.
As the boundary of such a neighborhood is a lens space $L(n,1)$, we have $\partial X=\Sigma \cup L(n,1)$. Now it
is easy to see that $X$ is $\Z_2$--acyclic, hence $\Sigma$ is $\Z_2$--homology cobordant to $-L(n,1)$. This implies
that $\Sigma$ has infinite order in $\Theta^3_{\Z_2}$ if $n \neq \pm 1$. Note that in the case that $n = \pm 1$, $X$ is
even a $\Z$--acyclic manifold and as $L(1,1)=S^3$ we obtain that the integral homology sphere $\Sigma$ is the
boundary of a $\Z$--acyclic manifold.
\end{example}

We close this section with the remark that the above computations for lens spaces can be generalized to
prove bounds for $m(\Sigma)$ and $\mbar(\Sigma)$ if $\Sigma$ is a Seifert fibred $\Z_2$--homology sphere, 
note that these spaces are double coverings of the 3--sphere branched along a Montesinos link. As the arguments
and calculations are very similar to the computations in the (special) case of lens spaces
we only state the result for one class of examples.

\begin{example} Let $n$ be a positive integer and consider the Seifert fibred space 
$\Sigma_n=M((12,5),(7,3),(6n+1,4n+1);-7)$, which is the double coverings of $S^3$
branched along the Montesinos knot $K={\mathfrak m}((12,5),(7,3),(6n+1,4n+1);-7)$. The obvious Seifert
surface for this knot is obtained from Seifert surfaces for the two--bridge knots 
$S(12,5),S(7,3)$ and $S(6n+1,4n+1)$ by adding three bands, one of them having $7$ half twists.
Using this and the results about the signature and the slice genus of a two--bridge knot which we obtained in 
this section one easily derives that $\sigma(K)=-2(n+3)$ and $g^*(K) \leq n+3$. Note that, as for every knot,
we have the inequality $2g^*(K) \geq |\sigma(K)|$, so our estimate for the slice genus happens to be sharp.

Now we can conclude that
$\mbar(\Sigma_n) \leq -\frac{1}{2}(n+3)$ . By Theorem~\ref{basic} this shows that $\Sigma_n$
has infinite order in $\Theta^3_{\Z_2}$. Also note that the order of the first homology of $\Sigma_n$
is linear in $n$, so $\Sigma_n \neq \Sigma_m$ if $m \neq n$, and we obtain an infinite family of Seifert
fibred spaces which have infinite order in~$\Theta^3_{\Z_2}$.
\end{example}

\section{A surgery formula}

As indicated in the introduction, a second main source of $\Z_2$--homology spheres is surgery on knots
with odd integral framings.
In this section, we shall see how one can obtain information on the invariants
$m(\Sigma)$ and $\mbar(\Sigma)$ if a $\Z_2$--homology sphere $\Sigma$ is described by integral surgery
on a knot. Combined with other methods for calculating these invariants, this provides a lower bound for 
the slice genus of such a knot. We also give examples of 3--manifolds which cannot be obtained by integral
surgery on a knot.

Suppose we are given a knot $K$ in the 3--sphere. Then surgery with odd framing on this knot yields a
$\Z_2$--homology sphere $\Sigma$.
First let us derive a relation between the Arf invariant of the knot, the framing and the Rokhlin invariant 
$R(\Sigma)$ of $\Sigma$. For this purpose recall that if $W$ is a simply connected
4--manifold whose boundary is a $\Z_2$--homology sphere and $F \subset W$ is a closed embedded characteristic surface, 
we have a well defined Arf invariant
$Arf(F) \in \Z_2$ which does only depend of the homology class of the surface and the Rokhlin invariant
of the boundary (as one can see by gluing with a simply connected spin manifold with boundary $-\partial W$,
note that the condition that $\partial W$ is a $\Z_2$--homology sphere implies that $F$ is still characteristic
in the resulting closed 4--manifold).

\begin{prop}\label{arf}
Let $K$ be a knot and assume that $n$ is some odd number. Let $\Sigma$ denote the $\Z_2$--homology sphere
which is the result of surgery along $K$ with framing $n$. Then
$n-\frac{n}{|n|}\equiv - R(\Sigma)\mod 8$ and 
\[
Arf(K) = \frac{1}{8}(n-\frac{n}{|n|} + R(\Sigma)) \mod 2.
\]
\end{prop}

\begin{proof}
Let $W$ denote the simply connected 4--manifold which is obtained from $D^4$ by adding a handle along $K$
with framing $n$. Then $\partial W=\Sigma$ is a $\Z_2$--homology sphere, in fact $\# H_1(\Sigma;\Z)=|n|$.
Furthermore $b_2(W)=1$. As the boundary of $W$ is a $\Z_2$--homology sphere, the
$\Z_2$--intersection form is non--degenerate and hence the fact that $b_2(W)=1$ implies that $W$ is not spin, i.e.
$w_2(W)$ is the non--zero element of $H^2(W;\Z_2)=\Z_2$. Pick a Seifert surface $F'$ for $K$ and let $F \subset W$
denote the surface which is obtained by gluing $F'$ with the core of the 2--handle. Then $F \cdot F=n$,
and therefore $sign(W)=\frac{n}{|n|}$. Let $[F] \in H_2(W;\Z_2)$ denote the homology class of $F$ and denote the
Poincar\'e duality map $H_2(W;\Z_2) \rightarrow H^2(W,\partial W;\Z_2)$ by $PD_W$. Then the fact that the
self--intersection of $F$ is odd implies that the image of $PD_W([F])$ under the restriction 
$\iota^* \colon H^2(W,\partial W;\Z_2) \rightarrow H^2(W;\Z_2)$ is the non--zero element, i.e. we have
\[
\iota^*PD_W([F])=w_2(W).
\]
Now pick a spin manifold $V$ with boundary $-\Sigma$ and consider the closed manifold $X=W \cup V$. The fact that
$\Sigma$ is a $\Z_2$--homology sphere implies that $H^2(X;\Z_2)=H^2(W;\Z_2) \oplus H^2(V;\Z_2)$. As $V$ is spin
the surface $F$ represents $w_2(X)$, i.e. $F \subset X$ is characteristic. Therefore we can conclude that
$F \cdot F \equiv sign(X) \mod 8$. By Novikov additivity, $sign(X)=sign(W) + sign(V)=\frac{n}{|n|} + sign(V)$.
By definition of the Rokhlin invariant we also have that $sign(V) \equiv -R(\Sigma) \mod 8$, and therefore we
obtain that $n \equiv \frac{n}{|n|}-R(\Sigma) \mod 8$. By \cite{FK}, we obtain
\[
Arf(F)=\frac{1}{8} (F \cdot F - sign(X)) = \frac{1}{8}(n - \frac{n}{|n|} + R(\Sigma)) \mod 2
\]
for the Arf invariant of the surface $F$. However it is also known~\cite{Ki} that $Arf(F)=Arf(K)$, 
and the proof of the proposition is complete.
\end{proof}

\begin{cor}\label{congruence}
Assume that a $\Z_2$--homology sphere $\Sigma$ is obtained by integral surgery on a knot.
Then
\[
\# H_1(\Sigma;\Z) - 1 \equiv \pm R(\Sigma) \mod 8.
\]
Here the sign is minus if the framing is positive, otherwise it is plus.
\end{cor}

\begin{proof}
Suppose that $\Sigma$ is obtained by surgery on a knot $K$ with framing $n$. First let us consider the case
that $n$ is positive (note that $n$ must be odd). Then $\# H_1(\Sigma;\Z)=n$, and by Proposition~\ref{arf},
we have that $\# H_1(\Sigma;\Z) - 1 \equiv - R(\Sigma) \mod 8$. If $n$ is negative, we obtain
$- \# H_1(\Sigma;\Z) +1 \equiv -R(\Sigma) \mod 8$, and multiplying this by minus one gives the desired
result.
\end{proof}

\begin{remark} If a lens space $L(p,q)$ can be obtained by surgery on a knot, then $q$ or $-q$ is 
a square modulo p by Proposition 1 in~\cite{FS1}. It is interesting that if $p$ is an odd prime,
this criterion is actually  eqivalent to the congruence of Corollary~\ref{congruence},
which can therefore be seen as a generalisation of the result for lens spaces in~\cite{FS1}.
In fact, assume that $q$ is a square mod $p$ (we can restrict ourselves to this case after possibly reversing
the orientation), i.e. $\binom{q}{p}=1$. By~\cite{HZ}, p. 137,
we have
\begin{equation}\label{hirza}
\binom{q}{p}+6p s(q,p) \equiv \frac{p+1}{2} \mod 4
\end{equation}
where $s(q,p)$ is a Dedekind sum. It is also known \cite{W} that $R(L(p,q))\equiv \pm 4p^2 s(q,p) \mod 8$
(the sign depending on orientation conventions),
note that, as $p$ is odd, the Rokhlin invariant is even and therefore $2p^2s(q,p)$ is an integer.
Multiplying equation~\eqref{hirza} by 2p we therefore obtain
\[
12 p^2 s(q,p) \equiv -4 p^2 s(q,p) \equiv p(p-1) \mod 8.
\]
But $p$ is an odd integer, hence $p^2 \equiv 1 \mod 8$, and we end up with
\[
\pm R(L(p,q)) \equiv 4p^2 s(q,p) \equiv p-1 = \# H_1(L(p,q);\Z)-1  \mod 8
\]
which is the prediction made by Corollary~\ref{congruence}.
A similar calculation shows that if the congruence of Corollary~\ref{congruence} holds for a
lens space $L(p,q)$, one of the Jacobi symbols $\binom{q}{p}$ and $\binom{-q}{p}$ must be one,
hence, as $p$ is a prime, $q$ or $-q$ is a square modulo $p$.
\end{remark}

In some cases Corollary~\ref{congruence} can be used to show that certain 3--manifolds are not
the result of integral surgery on a knot (although they can of course be obtained by integral
surgery along a link) or to determine the sign of the framing.

\begin{example}\label{niceex}
Let us consider an example of a connected sum where each summand is the result of surgery
on a knot but the sum is not. Let $\Sigma=L(3,1) \# L(7,1)$. Note that the first homology of
$\Sigma$ is cyclic of order 21. 
Clearly both summands can be obtained by surgery on knots, namely by $-3$ respectively $-7$
surgery on the unknot. However $R(\Sigma)=R(L(3,1)+R(L(7,1))=2+6=8$, and as $21 \not\equiv 0 \mod 8$ we can
conclude that $\Sigma$ is not the result of integral surgery on a knot. 
Note that, as $\Sigma$ can of course be obtained by integral surgery along the trivial two--component
link, its surgery number as defined in~\cite{Au} is two.
The same argument shows that
all manifolds of the form $L(8k+3,1)\#L(8k+7,1)$ have surgery number two.
\end{example}

\begin{example}
Suppose that $n \equiv 3 \mod 4$. The lens space $L(n,1)$ is the result of surgery with framing
$-n$ on the unknot. However Corollary~\ref{congruence} shows that this manifold cannot be obtained
by surgery on a knot with {\em positive} framing. In fact, $R(L(n,1))=n-1$,
$\# H_1(L(n,1);\Z)=n$, and if $n \equiv 3 \mod 4$, then $n-1 \not\equiv -(n-1) \mod 8$.
\end{example}

As $\pm 2$ is not a square modulo 5, Proposition 1 in \cite{FS1} implies that $L(5,2)$ is not the result
of surgery on a knot.
In~\cite{Au}, D. Auckly gave another argument for this fact 
which can be generalized to the following 

\begin{prop}\label{unknotting}
Assume that $K$ is a knot, $det(K) \neq \pm 1$, 
which can be unknotted using only one crossing change, and suppose that neither $2$ nor
$-2$ is a square modulo $|det(K)|$. Let $\Sigma$ denote the double covering of the 3--sphere branched along
$K$. Then $\Sigma$ is the result of rational surgery on a knot but cannot be obtained by integral surgery on
a knot. 
\end{prop}

\begin{proof}
Suppose for a moment that $\Sigma$ can be obtained as the result of integral surgery on some knot. The trace of
this surgery is a simply connected 4--manifold $X$ with second homology $H_2(X;\Z)=\Z$. If $n$ denotes the framing 
of the surgery the intersection form on $X$ is $n$ times the standard form on $\Z$. Using this
one easily sees that the linking form on $H_1(\Sigma;\Z)$ maps the generator of $H_1(\Sigma;\Z)$ to 
$\pm \frac{1}{|n|}$,
see for instance~\cite{Au} for a proof.
Note that $n=\pm \# H_1(\Sigma;\Z)=\pm det(K)$.

Now it has been pointed out in \cite{CL} that -- as a consequence of the fact that the unknotting number is one --
the linking form on $H_1(\Sigma;\Z)$ takes precisely the values
$k^2 \frac{2\epsilon}{|det(K)|}$ in $\Q / \Z$ for some sign $\epsilon \in \{1,-1\}$, in particular 
$\frac{2\epsilon}{|det(K)|}$ is in the image.
Hence there exist integers $a,b \in \Z$ such that 
\[
\pm a^2 \frac{1}{|det(K)|} =  b + \frac{2\epsilon}{|det(K|}.
\] 
Multiplying this by $|det(K)|$ leads to $\pm a^2 \equiv 2 \mod |det(K)|$ in contradiction to our assumption.
\end{proof}

\begin{example}
Let $K$ denote the Montesinos knot ${\mathfrak m}(1;(2,1),(3,2),(3,2))$ which is denoted by 
${\bf 8}_{21}$ in~\cite{Ka}.
The double covering $\Sigma$ of $S^3$ branched along $K$ 
is then a Seifert fibred space $M(1;(2,1),(3,2),(3,2))$, note that this 3--manifold is actually irreducible. 
One easily sees that $|det(K)|=15$ and $u(K)=1$ (the
reason being that the rational tangle $t(3,2)$ can be changed to the trivial tangle by one crossing change).
As $\pm 2$ is not a square modulo 15, 
Proposition~\ref{unknotting} applies and we obtain that $\Sigma$ cannot be obtained by integral surgery on a knot
although it is the result of rational Dehn surgery on some knot in $S^3$. Note that $\# H_1(\Sigma;\Z)=15$
and $\sigma(K)=2$, hence $\# H_1(\Sigma;\Z)-1 \equiv -R(\Sigma) \mod 8$, so Corollary~\ref{congruence} does not yield
this result.
\end{example}

If a $\Z_2$--homology sphere is the result of surgery on a link then the trace of this surgery will be a natural
choice for a 4--manifold bounded by $\Sigma$. 
Unfortunately this trace will in general not be spin, but one can always
do surgery on an embedded surface to obtain a spin manifold.

\begin{lemma}\label{ylike}
Suppose that $W$ is a simply connected 4--manifold whose boundary is a $\Z_2$--homology sphere 
and that $F \subset W \setminus \partial W$ is a closed characteristic surface 
having non--zero self--intersection number.
Let $\epsilon \in \{-1,+1\}$
denote the sign of $F \cdot F$. Then there exists a simply connected spin 4--manifold  
$X$ with $\partial X=\partial W$
such that 
\begin{align*}
b_2(X) &=b_2(W)+2(g(F)-1)+|F \cdot F + 8\epsilon Arf(F)| + 4 Arf(F)\\
\sigma(X) &=\sigma(W)-(F \cdot F+8\epsilon Arf(F))
\end{align*}
where we think of the Arf invariant $Arf(F) \in \Z_2$ as an element of $\{0,1\}$.
\end{lemma}

\begin{proof}
Clearly we can restrict ourselves to the case that the self--intersection number of $F$ is positive, the
case of negative self--intersection number then follows by reversing the orientation. 

Let us first consider the case
that the Arf invariant $Arf(F)$ is zero. Let $W'$ denote the manifold which is obtained from $W$ by attaching
$F \cdot F -1$ copies of  $\cpq$ and consider the surface $F' \subset W$ given by gluing
$F$ with the exceptional divisors. Then $F' \cdot F'=1$, and clearly $Arf(F')=Arf(F)=0$. As explained
in~\cite{Y} we can now construct a sphere $S \in W''$, where $W''$ is obtained from $W$ by attaching 
$g(F)$ copies of $S^2 \times S^2$ such that $S \cdot S = F' \cdot F'=1$ and such that $S$ is still characteristic.
Blowing down this sphere gives a spin 4--manifold $X$ as required.

In the case that $Arf(F)=1$, consider $W'=W \# \C P^2$ and the surface $F' \subset W'$
obtained by gluing $F$ with a torus representing 3 times the generator of $H_2(\C P^2;\Z)$.
Then $Arf(F')=0$ and the claim follows from what we just proved.
\end{proof}

\begin{prop}\label{estimate}
Suppose that a $\Z_2$--homology sphere is obtained by integral surgery with framing $n$ on a knot $K$. Let
$\epsilon=\frac{n}{|n|}$ denote the sign of the framing and let
\[
\mu=\frac{1}{8} (n-\epsilon) + R(\Sigma)) \! \! \mod 2 \in \Z_2 
\]
(note that, according to Proposition~\ref{arf}, the number in parentheses is actually a multiple of 8).
Finally let $h=\# H_1(\Sigma;\Z)$ and let $g^*=g^*(K)$ denote the slice genus of the knot $K$.
\begin{enumerate}
\item If $\mu=0$, then 
\[
\frac{4-5\epsilon}{4}(h-1) + 2g^* \geq \mbar(\Sigma) \geq m(\Sigma) \geq \frac{-4-5\epsilon}{4}
(h-1) - 2g^*.
\]
\item If $\mu=1$, then 
\[
\frac{4-5\epsilon}{4}(h+7) + 2g^* + 4 \geq \mbar(\Sigma) \geq m(\Sigma) \geq \frac{-4-5\epsilon}{4}(h+7)-2g^* -4.
\]
\end{enumerate}
\end{prop}

\begin{proof}
By assumption, $\Sigma$ is obtained by surgery on a knot $K$. The trace of this surgery is a
4--dimensional handlebody $W$ with boundary $\Sigma$ which is the result of attaching a
2--handle along $K$ to $D^4$. Consequently $b_2(W)=1$ and $\sigma(W)=\epsilon$. By
Proposition~\ref{arf}, $Arf(K)=\mu$, and clearly $|n|=h$, i.e. $n=\epsilon h$.

Now pick a surface $F \subset D^4$ with boundary $4\partial F=K$ 
such that $g(F)=g^*(K)$ and let $F'$ denote a Seifert surface
for $K$. Then the surfaces $G$ and $G'$ formed by gluing the core 
of the 2--handle with $F$ respectively
$F'$ clearly have the same homology class $\alpha$. As in the proof of Proposition~\ref{arf}, this homology class
is characteristic, and $\alpha \cdot \alpha=n$. It is known that the Arf invariant of the surface $G$ obtained by
gluing $F'$ and the core of the 2--handle is $Arf(K)$, and as $[G]=[G']$, the same is true for $G'$
Now we can apply
Lemma~\ref{ylike} to $W$ to obtain a spin 4--manifold $X$ with boundary $\Sigma$
and the claimed estimates follow from the
definitions of the invariants $m$ and $\mbar$.
\end{proof}

\begin{remark} A similar estimate can be derived in the more general case that
the $\Z_2$--homology sphere $\Sigma$ is obtained by integral surgery on a link which has a
characteristic knot (which one can always assume after sliding handles), such a characteristic knot
will again define a characteristic surface to which one can apply Lemma~\ref{ylike}.
\end{remark}

We are now ready to 
use Proposition~\ref{estimate}, combined with the information obtained from Proposition~\ref{arf},
to derive a lower bound for the slice genus of a knot on which integral surgery can be performed to obtain 
a given $\Z_2$--homology sphere $\Sigma$. Note that if $\Sigma$ is the result of surgery of a knot, then we can,
by Corollary~\ref{congruence}, choose the orientation of $\Sigma$ such that
$\# H_1(\Sigma;\Z) - 1 \equiv -R(\Sigma) \mod 8$, therefore we state our result only in this case.

\begin{thm}\label{knotestimate}
Assume that $\Sigma$ is a $\Z_2$--homology sphere and that $R(\Sigma) \not\equiv 4 \mod 8$
and $\# H_1(\Sigma;\Z) -1 \equiv -R(\Sigma) \mod 8$
Let 
\[
\mu=\frac{1}{8} (\# H_1(\Sigma;\Z) -1 + R(\Sigma)) \mod 2
\]
were we think of $\mu$ as an element of $\Z_2=\{0,1\}$.
If $\Sigma$ can be obtained by integral surgery on a knot $K$, then
\[
g^*(K) \geq \frac{1}{8}(\# H_1(\Sigma;\Z)-1 + 4 m(\Sigma))-\mu.
\]
\end{thm}

\begin{proof}
Let $n$ denote the framing of the knot. Of course $n=\pm \# H_1(\Sigma;\Z)$. Once we can show that $n$ must be
positive the claimed inequality follows from
Proposition~\ref{estimate}. So let us assume that $n < 0$. Then $n=-\# H_1(\Sigma;\Z)$, and by Proposition~\ref{arf},
we have
\[
-\# H_1(\Sigma) + 1 \equiv -R(\Sigma) \mod 8.
\]
However combining this with the assumption $\# H_1(\Sigma;\Z) -1 \equiv -R(\Sigma) \mod 8$
shows that $R(\Sigma) \equiv - R(\Sigma) \mod 8$. But this is only possible if $R(\Sigma) \equiv 0 \mod 4$,
which contradicts our assumptions. 
\end{proof}

\begin{example}\label{series}
Assume that $k$ is an even number and consider the lens space $L(16k+7,7k+3)$. Note that this lens space
is -- up to orientation -- $L(4mn-1,4n^2)$ with $n=2k+1$, $m=2$. As $n$ and $m$ are coprime, $L(16k+7,7k+3)$
can be obtained by surgery on a knot by Theorem 1 in~\cite{FS1}. 
An admissible continuous
fractions decomposition is given by
\[
\frac{16k+7}{7k+3}=[2,4,-1,-2,1,k,-1].
\]
Using this one can easily derive that $R(L(16k+7,7k+3))=2$ and $m(L(16k+7,7k+3))\geq -1.5$. Hence
Theorem~\ref{knotestimate} implies that the slice genus of any knot $K$ on which integral surgery can be performed
to obtain
$L(16k+7,7k+3)$ must be at least $2k-1$.
Of course the slice genus of a knot is at most its genus,
and for every knot the crossing number is at least two times the genus
(otherwise Seifert's algorithm would produce a Seifert surface of smaller genus), so we also obtain
a lower bound for the crossing number.
\end{example}

The most obvious infinite family of lens spaces which can be obtained by integral
surgery on a knot is given by the lens spaces $L(n,1)$. All these spaces can be obtained by
surgery on a single knot, namely the unknot (but of course with different framings). Example~\ref{series}
shows that this situation is not typical and that one actually has to use infinitely many knots (even
infinitely many concordance classes of knots) to obtain all the
lens spaces which are the result of integral surgery on a knot, i.e. we have the

\begin{cor}
For every natural number $N > 0$
there exists a lens space $L(p,q)$ which can be obtained by integral surgery on a knot
such that every such knot has slice genus at least $N$. 
\end{cor}

\end{document}